\documentclass[11pt]{article}
\usepackage{shapepar,bm}
\usepackage[dvips,arrow,matrix,ps,color,line]{xy}
\usepackage{theorem,amsfonts,amssymb,amscd}
\usepackage{geometry}
\usepackage{graphicx}
\usepackage{graphics}
\usepackage{amssymb}

\usepackage{amssymb, graphics}

\newcommand{\mathsym}[1]{{}}

\usepackage[usenames]{color}
\definecolor{MyLightMagenta}{cmyk}{0.1,0.8,0,0.1}
\definecolor{MyDarkBlue}{rgb}{0.1,0,0.3}
\makeindex

\def\NN{\mathbb N}

\def\ZZ{\mathbb Z}
\def\CC{\mathbb C}

\def\PP{\mathbb P}

\def\TT{{\bf T}}

\def\ttp{{\tt p}}

\def\ovD{\overline{D}}

\def\ikfornu{{\nu^{i_1}\wedge\ldots\wedge\nu^{i_k}}}
\def\kfornu{{\nu^{1}\wedge\ldots\wedge\nu^{k}}}

\def\ep{{\epsilon}}

\def\ev{{\rm ev}}

\def\Bw{{\bigwedge}}

\def\w2M{\bigwedge^2M}

\def\wM{\bigwedge M}

\def\w{\wedge }
\def\bw{\bigwedge }

\protect\def\wkM{{\bigwedge^kM}} 
\protect 
\protect
\protect

\protect
\protect

\def\sra{\rightarrow}
\def\lra{\longrightarrow}

\def\proof{\noindent{\bf Proof.}\,\,}
\def\qed{{\hfill\vrule height4pt width4pt depth0pt}\medskip}
\def\be{\begin{equation}}
\def\ee{\end{equation}}
\def\bclm{\begin{claim}}
\def\eclm{\end{claim}}
\def\beqn{\begin{eqnarray}}
\def\eeqn{\end{eqnarray}}
\def\beqn*{\begin{eqnarray*}}
\def\eeqn*{\end{eqnarray*}}

\def\kformX{{X^1\w\ldots\w X^{k}}}

\def\kformep{{\epsilon^1\wedge\ldots\wedge\epsilon^k}}

\def\ikformep{\epsilon^{i_1}\wedge\ldots\wedge\epsilon^{i_k}}

\theoremstyle{change}
\theorembodyfont{\rmfamily}

\newtheorem{claim}{}[section]

\global
\global

\def\no@breaks#1{{\def\\{ \ignorespaces}#1}}    % disallow explicit line breaks

  \makeatletter
\def\cleardoublepage{\clearpage\if@twoside \ifodd\c@page\else
\hbox{} \thispagestyle{empty}
\newpage
\if@twocolumn\hbox{}\newpage\fi\fi\fi} \makeatother

\usepackage{eso-pic,graphicx}
\makeatletter
\newcommand\BackgroundPicture[2]{%
  \setlength{\unitlength}{1pt}%
  default \put(0,\strip@pt\paperheight){%
  \parbox[t][\paperheight]{\paperwidth}{%
    \vfill
     \centering \includegraphics[angle=#2, width=15cm, height=15cm,  bb=0 0 150 150]{#1}
    \vfill
}}} %
\makeatother
\usepackage{geometry}                % See geometry.pdf to learn the layout options. There are lots.
\usepackage{graphicx}
\usepackage{amssymb}
\usepackage{epstopdf}

\title{Schubert Calculus on a Grassmann Algebra\thanks{\noindent Work partially sponsored by PRIN ``Geometria sulle Variet\`a Algebriche" (Coordinatore 
A.~Verra), INDAM-GNSAGA and ScuDo, Politecnico di Torino.}}

\author{Letterio Gatto, Ta\'\i se Santiago }

\date{}                                           % Activate to display a given date or no date

\begin{document}

\maketitle

\abstract{\small
\noindent 
The ({\em classical},  {\em small quantum}, {\em equivariant}) cohomology ring of the grassmannian $G(k,n)$ is generated by certain derivations operating on an exterior algebra of a free module  of rank $n$ ({\em Schubert Calculus on a Grassmann Algebra)}. Our main result gives, in a unified way, a presentation of all such cohomology rings in terms of generators and relations. It also provides, by results of  Laksov and Thorup (\cite{LakTh} and \cite{LakTh1}), a presentation of the universal splitting algebra of a monic polynomial of degree $n$ into the product of two monic polynomials, one of degree $k$.
 }

\section{Introduction}  In the paper~\cite{Gat1}  one shows that the cohomology ring of the complex grassmannian $G(k,n)$, parametrizing $k$-dimensional subspaces of $\CC^n$, can be realized as a commutative ring of endomorphism of the $k^{th}$ exterior power of a free $\ZZ$-module $M$ of rank $n$. Such a result was achieved  by studying a natural {\em Hasse-Schmidt derivation} on  the exterior algebra of $M$; Laksov and Thorup~(\cite{LakTh} and~\cite{LakTh1})  generalized it to the more interesting situation regarding the cohomology of Grassmann bundles.
 Their point of view is quite different, as it is based on the fact that the $k^{th}$-exterior power of a free $A$-module of rank $n$ can be endowed with a  natural module structure over the ring of symmetric polynomials (with $A$-coefficients):  this   leads to a beautiful  and natural description of the cohomology of $G(k,E)$, the  Grassmann bundle of $k$-dimensional subspaces in the fibers of a vector bundle $E$,  in terms of the {\em universal splitting algebra} of a certain monic polynomial $\ttp$ (encoding the Chern classes of $E$) into the product of two monic polynomials, one of degree $k$ (Cf. remark~\ref{rmklakth}). 
 
The main goal of this paper is to generalize~\cite{Gat1} via a translation of Laksov and Thorup's formalism into the language of derivations. 
 A {\em derivation}   on $\wM$,  the exterior algebra of a module $M$ over a commutative ring with unit,  is a sequence $D:=(D_0,D_1,\ldots)$ of endomorphisms, such that the $h^{th}$ order Leibniz's rule:
\be
D_h(\alpha\w\beta)=\sum_{\scriptsize \matrix{{h_1+h_2=h}\cr{ h_i\geq 0}}}D_{h_1}\alpha\w D_{h_2}\beta,\label{eq:hleibrul}
\ee
holds for each $h\geq 0$ and each $\alpha,\beta\in\wM$ (see~\ref{ref3.2}).
In~\cite{tesi}, any such a  derivation is  called a {\em Schubert Calculus on a Grassmann Algebra}. The terminology is motivated by the fact that if one takes $M$ to be a finite free module over a graded commutative $\ZZ$-algebra  of characteristic $0$, there is a {\em canonical derivation} on $\wM$ (generalizing that studied in~\cite{Gat1}; see Section~\ref{sect2}) describing, within a unified framework, different kind of cohomology theories on  complex grassmannian varieties,  such as, e.g., the classical, the small quantum or the equivariant one. Working  on the exterior algebra, instead of on a single exterior power,  many formal manipulations get easier: as an example we offer Theorem~\ref{generalpres}, the main result  of this paper, that consists in a simple formula {\em giving, in a unified way,  the presentation of the classical, small quantum and equivariant cohomology ring of the complex  grassmannian} $G(k,n)$. In fact,  the (classical, small quantum, equivariant) cohomology ring of  all the grassmannians $G(k,n)$, $1\leq k\leq n$, are quotient of a same commutative ring of endomorphisms of  the exterior algebra of a free module of rank $n$ (see Sect.~\ref{ref4.7}). As the latter is generated by derivations, the (classical, small quantum, equivariant)  Schubert calculus on $G(k,n)$  can be reduced to that, much easier,  on $G(1,n)=\PP^{n-1}$ (as  in~\cite{Gat1}; see also~\cite{Gat2}).   Our best application of such a philosophy regards an elementary description, as in~\cite{GatSant2} (see also~\cite{tesi}), of the equivariant Schubert calculus  on a grassmannian acted on by a torus with isolated fixed locus,   recovering, in particular,  the case studied  in~\cite{KT} (see also~\cite{formalism}).

\medskip
{\bf Acknowledgment.}  The first author wants to thank the warm ospitality of the STID of Menton, Universit\'e de Nice, Sophia-Antipolis, notably that of  its chairmain, Guy Choisnet, where most part of this paper, originated from~\cite{Gat1} and~\cite{tesi}, has been written. The current exposition  has been deeply influenced by the work of D.~Laksov and A.~Thorup on related subjects (\cite{LakTh}, \cite{LakTh1}, \cite{formalism}) and by many  conversations the authors had with the former, to whom they want to address a warm feeling of gratitude. We also thank I.~Vainsencher for some key suggestions as well as the Referee for his valuable and (especially) patient remarks.

%\section{Preliminaries and Notation} 
%\claim{}\label{notpih} Let ${\cal I}^k=\{I=(i_1,\ldots,i_k)\in\NN^k\,|\,1\leq i_1<\ldots<i_k\}$ (as in ~\cite{anderson}, \S 5, Section~1) and~\cite{Gat1}). The {\em weight} of $I\in{\cal I}^k$ is $wt(I)=\sum_{j=1}^k(i_j-j)$. It coincides with the weight of the associated  partition $(i_k-k, i_{k-1}-(k-1),\ldots,i_1-1)$. 
%If $I\in{\cal I}^k$, let
%\[
% {\cal P}(I):=\{H=(h_1,\ldots,h_k)
%\in\NN^k\,\,|\,\,i_1\leq i_1+h_1<i_2\leq \ldots\leq i_{k-1}+h_{k-1}<i_k\}
%\]
%and ${\cal P}(I,h)=\{(h_1,\ldots,h_k)\in{\cal P}(I)\,\,|\,\, \sum_{i=1}^kh_i=h\}$.
% \claim{} \label{prel2} If $A$ is a commutative ring with unit, we denote by $A[\TT]$  the ring of polynomials in infinitely many indeterminates  $\TT:=(T_1,T_2,\ldots)$ and  by $A[\TT_k]$ the subring $A[T_1,\ldots,T_k]$. 
%As usual, $\Delta_I(\TT)\in A[\TT]$ will stand for the Schur polynomial  $\Delta_{(i_1,\ldots,i_k)}(\TT)=\det(T_{i_j-i})$ (setting $T_0=1$ and $T_j=0$, if $j<0$). Notice that $\Delta_{(23\ldots k+1)}(\TT)$ lands in fact in $A[\TT_k]$. Expanding the determinant along the last column, one sees  that $\Delta_I(\TT)$ belongs to the ideal $(T_{i_k-1},\ldots, T_{i_k-k})$ of $A[\TT]$.

%
% If  $A:=\bigoplus_{i\geq 0}A_i$ is a graded ring and  $a\in A_l$,  the {\em degree}  of the monomial $aT_{i_1}^{m_1}\ldots T_{i_j}^{m_j}$  is {defined to be} $l+m_1i_1+\ldots+m_ji_j$. Then  $A[\TT]$ is itself a graded ring $\bigoplus_{h\geq 0}A[\TT]_h$, where $A[\TT]_h$ is the submodule of all elements of $A[\TT]$ of degree $h$.
 
\section{Derivations on Exterior Algebras} \label{sect1}
\claim{} \label{prel3} Let $M$ be an $A$-module, $A[[t]]$ be the ring of  formal power series in an indeterminate $t$ over $A$ and  $\wM[[t]]:=(\wM)[[t]]$ be the $A[[t]]$-module of formal power series with coefficients in  $\wM=\bigoplus_{k\geq 0}\wkM$, the exterior algebra of $M$. The { former} gets a structure of $A[[t]]$-algebra by setting $
\sum_{i\geq 0}\alpha_it^i\wedge \sum_{j\geq 0}\beta_jt^j=\sum_{h\geq 0}\sum_{i+j=h}(\alpha_i\w\beta_j)t^h.
$

\claim{} \label{ref3.2} An $A$-module homomorphism $D_t:\wM\sra \wM[[t]]$ is said to be a {\em derivation} on $\wM$ if it is an $A$-algebra homomorphism, i.e. if  for each $\alpha,\beta\in\wM$:
\be
D_t(\alpha\w\beta)=D_t\alpha\w D_t\beta.\label{eq:fundeq}
\ee
The algebra homomorphism $D_t$ can be written as  a formal power series   $\sum_{i\geq 0}D_it^i$, with coefficients in the $A$-algebra $End_A(\wM)$. Denote by $D$ the sequence $(D_0,D_1,\ldots)$ of the coefficients of $D_t$. Equation~(\ref{eq:fundeq}) implies that for each $h\geq 0$, the $A$-endomorphism $D_h$ of $\wM$ satisfies the $h^{th}$-order Leibniz rule~(\ref{eq:hleibrul}),
got by expanding both sides of~(\ref{eq:fundeq}) and equating the coefficients of $t^h$ occurring on both sides. 
\claim{} Let $\jmath:Hom_A(\wM,\wM[[t]])\sra End_A(\wM[[t]])$ be the natural map sending any $\Psi_t=\sum_{i\geq 0}\psi_it^i\in Hom_A(\wM,\wM[[t]])$ to the endomorphism $\jmath(\Psi)$ of $\wM[[t]]$, defined, on each  $\sum_{i\geq 0}\alpha_it^i\in\wM[[t]]$, as:
\[
\jmath(\Psi)(\sum_{i\geq 0}\alpha_it^i)=\sum_{i\geq 0}\Psi(\alpha_i)\cdot t^i=\sum_{h\geq 0}\big(\sum_{i+j=h}\psi_i (\alpha_j)\big)t^h.
\]
If $D_t$ is a derivation,  then $\jmath({D}_t)$ is itself an $A[[t]]$-algebra endomorphism of $\wM[[t]]$. In fact it is obviously an $A[[t]]$-module endomorphism and, moreover:
\begin{eqnarray}
\jmath(D_t)\Big(\sum_{i\geq 0}\alpha_it^i\wedge\sum_{j\geq 0}\beta_jt^j\Big)&=&\jmath(D_t)\sum_{h\geq 0}\Big(\sum_{i+j=h}\alpha_i\w\beta_j\Big)t^h=\nonumber\\=\sum_{h\geq 0}\Big(\sum_{i+j=h}D_t(\alpha_i\w\beta_j)\Big)t^h&=&\sum_{h\geq 0}\Big(\sum_{i+j=h}D_t\alpha_i\w D_t\beta_j\Big)t^h=\nonumber\\=\sum_{i\geq 0}D_t\alpha_i\cdot t^i\wedge \sum_{j\geq 0}D_t\beta_j\cdot t^j&=&\jmath(D_t)\sum_{i\geq 0}\alpha_i\cdot t^i\w\jmath(D_t)\sum_{j\geq 0}\beta_j\cdot t^j.\label{eq:dcapaom}
\end{eqnarray}

\claim{} For each pair $D_t, D_t'\in  Hom_A(\wM,\wM[[t]])$, define a product $D_t*D'_t$  through the equality: 
$
(D_t*D'_t)\alpha=\jmath(D_t)(D'_t\alpha).
$
Clearly $\jmath(D_t)\alpha=D_t\alpha$  for each $\alpha\in\wM$ and
\begin{eqnarray}
(D_t*D'_t)(\alpha)&=&\sum_{h\geq 0}\big(\sum_{i+j=h}D_i(D'_j\alpha)\big)t^h=\jmath(D_t)(\sum_{j\geq 0}D'_j\alpha\cdot t^j)=\nonumber\\&=&\jmath(D_t)(D'_t\alpha)=(\jmath(D_t)\circ\jmath({D}'_t))\alpha.\label{eq:ddcap}
\end{eqnarray}

\claim{}\label{dt*dt'} The product $D_t*D'_t$ of two derivations on $\wM$  is a derivation on $\wM$.
Indeed, using~(\ref{eq:dcapaom}) and~(\ref{eq:ddcap}):
\begin{eqnarray*}
(D_t*D'_t)(\alpha\w\beta)&=&\jmath(D_t)(D'_t(\alpha\w\beta))=\jmath(D_t)(D'_t\alpha\w D'_t\beta)=\\&=&\jmath(D_t)(D'_t\alpha)\w \jmath(D_t)(D'_t\beta)=(D_t*D'_t)\alpha\w (D_t*D'_t)\beta,
\end{eqnarray*}
as desired. Let now $D^{(1)}=(D_i^{(1)})_{i\geq 0}$ be any (possibly finite) sequence of endomorphisms of $M$ and, for each $m\in M$, let $D_t^{(1)}(m)=\sum_{i\geq 0}D_i^{(1)}(m)t^i$. Then
$
D_t^{(1)}:M\sra M[[t]]
$
 is an $A$-module homomorphism. 
 
\bclm{\bf Proposition.} \label{extder} {\em There exists a unique derivation  $D_t:\wM\sra\wM[[t]]$ such that ${D_t}_{|_M}=D_t^{(1)}$ (or, equivalently, ${D_i}_{|_M}=D_i^{(1)}$).}\eclm

\proof 
For each $k\geq 1$, consider the $A$-multilinear map
$
M^{\otimes k}\sra(\wkM)[[t]]
$
defined by
$
 m_{i_1}\otimes\ldots\otimes m_{i_k}\mapsto D_t^{(1)}m_{i_1}\wedge\ldots\wedge D_t^{(1)}m_{i_k}
$, 
which is clearly alternating. By the universal property of exterior powers,  it factors through a unique $A$-module homomorphism  $\wkM\sra (\wkM)[[t]]$, given   by
$
D_t^{(k)}(m_{i_1}\wedge\ldots\wedge m_{i_k})=D_t^{(1)} m_{i_1}\wedge\ldots\wedge D_t^{(1)} m_{i_k}
$  on  monomials. Let $D_t\alpha=D_t^{(k)}\alpha
$ for all $\alpha\in\wkM$ and  all $k\geq 0$.
It follows that if $\alpha\in \bigwedge^{k_1}M$ and $\beta\in \bigwedge^{k_2}M$, equation~(\ref{eq:fundeq}) holds
by definition of $D_t$ and the fact that  $\alpha\wedge\beta$ is a finite $A$-linear combination of elements of the form
\begin{center}
$
\{m_{i_1}\wedge\ldots\wedge m_{i_{k_1}}\wedge m_{i_{k_1+1}}\wedge\ldots\wedge m_{i_{k_1+k_2}}; \quad 1\leq  i_1<\ldots<i_{k_1+k_2}\}.
$
\end{center}

\noindent
Since any  element of $\wM$ is a finite sum of homogeneous ones,  equation~(\ref{eq:fundeq}) holds for any arbitrary pair as well. The unicity part is straightforward: were $D'_t$ another extension of $D_t^{(1)}$, one would have
$
D_t'(  m_{i_1}\wedge\ldots\wedge m_{i_k})= D_t^{(1)} m_{i_1}\wedge\ldots\wedge D_t^{(1)} m_{i_k}=D_t( m_{i_1}\wedge\ldots\wedge m_{i_k}),
$
 for each $  m_{i_1}\wedge\ldots\wedge m_{i_k}$ and each $k\geq 1$.  Hence $D_t'=D_t$. 
 \qed

\claim{} Let ${\cal S}_t(\wM)$ be the set of all derivations $D_t:=\sum_{i\geq 0} D_it^i$ such that ${D_i}_{|_M}\in End_A(M)$ (i.e. the submodule $M$ of $\wM$ is $D_i${-stable}) and ${D_0}_{|_M}$  is an isomorphism. Hence $D_0:\wM\sra \wM$ is an isomorphism too.

\bclm{\bf Proposition.} {\em The pair  $({\cal S}_t(\wM),*)$ is a group.}\eclm

\proof
By \ref{dt*dt'}, ${\cal S}_t(\wM)$ is closed under $*$. By its very definition, $*$ is associative. The map ${\bf 1}:\wM\sra(\wM)[[t]]$, sending any $\alpha\in \wM$ to itself, thought of as a constant formal power series, is the $*$-neutral element. Thinking to $D_t$  as  a formal power series with coefficients in $End_A(\wM)$, the formal inverse ${D^{-1}_t}$  of $D_t$ (existing because of the invertibility of $D_0$) is a derivation as well. In fact
\begin{eqnarray*}
{ D^{-1}_t}(\alpha\w\beta)&=&\jmath({{D}}^{-1}_t)\big((D_t*{  D^{-1}_t})\alpha\w(D_t*{  D^{-1}_t})\beta)=\\&=&\jmath({{D}}^{-1}_t)(\jmath(D_t){  D^{-1}_t}\alpha\w\jmath(D_t){  D^{-1}_t}\beta)=\\&=&(\jmath({D}^{-1}_t)\circ\jmath(D_t))({  D^{-1}_t}\alpha\w{  D^{-1}_t}\beta)={  D^{-1}_t}\alpha\w {  D^{-1}_t}\beta,
\end{eqnarray*}
{  since} ${  D^{-1}_t}*D_t=D_t*{  D^{-1}_t}={\bf 1}$.
\qed

\claim{}\label{convx3.8} We fix another piece of notation. Let $A[\TT]$ be the polynomial ring in infinitely many indeterminates $\TT=(T_1,T_2,\ldots)$.  For each $k$-tuple $I:=(i_1,\ldots, i_k)$ of positive integers, we denote by $\Delta_I(\TT):=\Delta_{(i_1,\ldots,i_k)}(\TT)$ the {\em Schur polynomial} $\det[(T_{i_j-i})_{1\leq i,j\leq k}]\in A[\TT]$ (setting $T_0=1$ and  $T_j=0$, if $j<0$). By expanding $\Delta_I(\TT)$ along the last column,  one sees that $\Delta_I(\TT)$ belongs to the ideal  $(T_{i_k-1},\ldots, T_{i_k-k})$ of $A[\TT]$. In particular $\Delta_{(2,3,\ldots,h+1)}(\TT)\in (T_1,\ldots, T_h)$. If $D:=(D_0,D_1,\ldots, )$ is the sequence of coefficients of some $D_t\in{\cal S}_t(\wM)$ such that $D_0=id_{\wM}$, one defines $\Delta_I(D)$ to be the evaluation of $\Delta_I(\TT)$ at $D$ (via the substitution $T_i\mapsto D_i$).

\claim{}\label{conv3.9} For each $i\geq 0$, define $\ovD_i\in End_A(\wM)$ via the equality
${  D^{-1}_t}=\sum_{i\geq 0}(-1)^i{\ovD}_it^i
$.
By equating the coefficients of the same power of $t$ on both sides of the equation \linebreak $D_t*{  D^{-1}_t}=1$, one gets $\ovD_0=D_0^{-1}$, while, for each $h\geq 1$:
\be
\ovD_h-\ovD_{h-1}D_1+\ldots+(-1)^hD_h=0, \label{eq:invDt1}
\ee
so that, e.g., $\ovD_1=D_1$, $\ovD_2=D_1^2-D_2$. In general, one has (see~\cite{Fu1}, Appendix A):
\be
\ovD_h=\Delta_{(2,3,\ldots,h+1)}(D).\label{eq:invDt2}
\ee

 \bclm{\bf Proposition (Integration by parts).} {\em  Let $D_t\in{\cal S}_t(\wM)$. Then:
\begin{small}
\be
D_h\alpha\wedge\beta=\sum_{i\geq 0}(-1)^iD_{h-i}(\alpha\wedge \ovD_i\beta)=D_h\alpha\wedge\beta-D_{h-1}\alpha\wedge \ovD_1\beta+\ldots+(-1)^iD_0\alpha\wedge\ovD_h\beta.\label{eq:intpart2}
\ee
\end{small}}
\eclm

\proof
One expands both sides of the equality
$
\jmath(D_t)(\alpha\wedge {  D^{-1}_t}\beta)=D_t\alpha\wedge\beta
$, 
and then compares the coefficients of $t^h$ occurring on each  side. \qed

{  \claim{{  \bf Example.}}} One has
$
D_1\alpha\wedge D_0\beta=D_1(\alpha\wedge\beta)-D_0\alpha\wedge\ovD_1\beta$ and:
\[
D_2\alpha\wedge D_0\beta=D_2(\alpha\wedge\beta)-D_1(\alpha\wedge\ovD_1\beta)+D_0\alpha\wedge\ovD_2\beta.
\]

\section{Schubert Calculus on a Grassmann Algebra}\label{sect2}

\claim{}\label{recalling01}  From now on,  $A$ will  be assumed to be any graded ring   $\bigoplus_{i\geq 0}A_i$ such that $A_0=\ZZ$. Let  $X$ be an indeterminate over $A$, $M:=XA[X]$ and  $M(\ttp):=M/\ttp M$, where $\ttp$ is either the $0$ polynomial or  a monic polynomial $X^n-e_1X^{n-1}+\ldots +(-1)^ne_n\in A[X]$ such that $e_i\in A_i$.   Then $M(\ttp)$ is a free $A$-module generated by  ${\bm\ep}=(\ep^i)_{1\leq i\leq n}$, where $n$ is either $\deg(\ttp)$ if $\ttp\neq 0$, or $\infty$ if  $\ttp= 0$.%%%%%%%%%%%%%%%%%%%%%%%%%%%%%%%%%%%%%%%%%%%%%%%%%%%

%%%%%%%%%%%%%%%%%%%%%%%%%%%%%%%%%%%%%%%%%%%%%%%%%%%
\claim{}\label{not23} Let ${\cal I}^k=\{I=(i_1,\ldots,i_k)\in\NN^k\,|\,1\leq i_1<\ldots<i_k\}$ (as in ~\cite{anderson}, \S 5, Section~1, and~\cite{Gat1}). The {\em weight} of $I\in{\cal I}^k$ is $wt(I)=\sum_{j=1}^k(i_j-j)$. It coincides with the weight of the associated  partition $(i_k-k, i_{k-1}-(k-1),\ldots,i_1-1)$.  If $I:=(i_1,\ldots,i_k)\in{\cal I}^k$, let $\w^I{\bm\ep}$ denote $\ikformep$. Each exterior power $\wkM(\ttp)$ is a free $A$-module with basis  $\Bw^k{\bm\ep}:=\{\w^I\ep:I\in{\cal I}^k_n\}$. If $a\in A_h$, the {\em weight} of $a\cdot\w^I{\bm\ep}$ is,  by definition,  ${  h}+wt(I)$. 
Set $(\bw^kM({\tt p}))_w=\bigoplus_{0\leq h\leq w}\big(\bigoplus_{wt(I)=h}A_{w-h}\cdot\w^I{\bm \ep}\big)
$.
Then  $\wkM(\ttp)=\bigoplus_{w\geq 0}(\bw^k M({\tt p}))_w$,  a graded $A$-module via {\em weight}.
\claim{}  \label{rec24} By~Proposition~\ref{extder} there is a unique sequence $D:=(D_0,D_1,\ldots)$ of $A$-endomor- phisms of $\wM(\ttp)$ such that i) (the $h^{th}$-order) {\em Leibniz's rule}~(\ref{eq:hleibrul}) 
holds for each $h\geq 0$ and each $\alpha,\beta\in\wkM(\ttp)$ and ii) the {\em initial conditions}  $D_h\ep^i=\ep^{i+h}$ are satisfied,  for each $h\geq 0$ and each   $ i\geq 1$. Notice that $D_i\circ D_j=D_j\circ D_i$ in $End_A(\wM(\ttp))$, as a simple induction shows.

\bclm{\bf Proposition.} \label{prepierprop}{\em The following formula holds:
\be
D_h\ikformep=\sum \ep^{i_1+h_1}\w\ldots \w \ep^{i_k+h_k},\label{eq:prepier}
\ee
the sum over all $h$-tuples  $(h_i)_{1\leq i\leq k}$ of non negative integers such that $h_1+\ldots+h_k=h$.
}
\eclm
\proof See~\cite{LakTh} or, since equality~(\ref{eq:prepier}) is defined over the integers, use the same inductive proof as in~\cite{Gat1}, Proposition~2.3.\qed

{\claim{\bf Example.}} When expanding $D_h\ep^{i_1}\w\ldots\w \ep^{i_k}$,  cancellations may occur on the right hand side of~(\ref{eq:prepier}), due to the $\ZZ_2$-symmetry of the $\wedge$-product. For instance:
\[
D_2(\ep^1\w \ep^2)=\ep^3\w \ep^2+\ep^2\w \ep^3+\ep^1\w \ep^4=\ep^1\w \ep^4.
\]
The surviving summands are predicted by {\em Pieri's formula} for $D_h$, a rule to speed up computations of ``derivatives" of $k$-vectors. 
\begin{claim}{\bf Theorem (Pieri's formula).} \label{pieri3.8} {\em {\em Pieri's formula} holds:
\be
D_h(\ep^{i_1}\w\ldots\w \ep^{i_k})=\sum_{{(h_i)\in P(I,h)}}
 \ep^{i_1+h_1}\w\ldots \w \ep^{i_k+h_k}, \label{eq:pieruno}
\ee
where, if $I=(i_1,\ldots,i_k)\in{\cal I}^k$, we denote by ${\cal P}(I,h)$  the set of all  $k$-tuples of non negative integers $(h_1,\ldots,h_k)$ such that 
$
 i_1+h_1<i_2\leq i_2+h_2<\ldots<i_{k-1}\leq i_k
$
and $h_1+\ldots+h_k=h$.
}

\end{claim}
\proof See~\cite{LakTh} or, since formula~(\ref{eq:pieruno}) is defined over the integers, use the same proof as in~\cite{Gat1}, Theorem~2.4.

%\claim{\bf Remark.} The reason why~(\ref{eq:pieruno}) is referred to as Pieri's formula, is that it coincides with the {\em combinatorial Pieri's formula} involving Young diagrams. See~\cite{Gat1}, Remark 2.8.

%

\claim{} \label{ref4.7} Let $A$ be as in~\ref{recalling01} and $A[\TT]$ be as in~\ref{conv3.9}. If $a\in A_l$,  the {\em degree}  of the monomial $aT_{i_1}^{m_1}\ldots T_{i_j}^{m_j}$  is {  defined to be} $l+m_1i_1+\ldots+m_ji_j$. Then  $A[\TT]$ is itself a graded ring $\bigoplus_{h\geq 0}A[\TT]_h$, where $A[\TT]_h$ is the submodule of all elements of $A[\TT]$ of degree $h$. There is a natural evaluation map,  $\ev_D:A[\TT]\rightarrow End_A(\wM(\ttp))$,  sending $P\in A[\TT]$ to $P(D)$ (got by ``substituting"  $T_i\mapsto D_i$ into $P$). We denote by  ${\cal A}^*(\wM(\ttp))$ the image  of $\ev_D$ in $End_A(\wM(\ttp))$ and by ${\cal A}^*(\wkM(\ttp))$ the image of the natural restriction map  
\[
\rho_k:{\cal A}^*(\wM(\ttp))\rightarrow End_A(\wkM(\ttp)),
\] 
given by $P(D)\mapsto (D)_{|_{\wkM(\ttp)}}$. 
Pieri's formula implies {\em Giambelli's  formula}, a special case of the general determinantal formula stated in~\cite{LakTh},  Main Theorem, which reads, in this case, as:
\be
\ep^{i_1}\w\ldots\w \ep^{i_k}=\Delta_{(i_1\ldots i_k)}(D)\cdot \ep^1\w\ldots \w \ep^{k},\label{eq:giambform}
\ee
where, as in~\ref{conv3.9},  
$
\Delta_{(i_1\ldots i_k)}(D)=\ev_D(\Delta_{(i_1\ldots i_k)}(\TT))$.
We have hence {  shown} that:
\begin{claim}{\bf Theorem.}\label{thm1rec} {\em The natural evaluation map $\ev_\kformep:{\cal A}^*(\wM(\ttp))\rightarrow\wkM(\ttp)$, mapping $P(D)\mapsto P(D)\kformep$ is surjective.}\qed
\end{claim}
\claim{}\label{ref4.10} It follows that $\ker(\rho_k)=\ker(\ev_\kformep)$ and then:
\be
{\cal A}^*(\wkM(\ttp))={{\cal A}^*(\wM(\ttp))\over\ker(\ev_\kformep)}.
\ee
The induced map $\Pi_k:{\cal A}^*(\wkM(\ttp))\sra \wkM(\ttp)$, defined by 
\[
P(D)+\ker\ev_\kformep\mapsto P(D)\kformep,
\] 
we call the
{\em Poincar\'e isomorphism}.

\claim{\bf Remark.} \label{prel211}  Let  ${\cal I}^k_n=\{I\in{\cal I}^k\,\,|\,\, i_k\leq n\}$.
 A routine check  shows that if $I=(i_1,\ldots, i_k)\in {\cal I}^k_n$ and $H\in{\cal P}(I,h)$ then $I+H:=(i_1+h_1,\ldots,i_k+h_k)\in {\cal I}^k$. Denote by ${\cal I}^{k,w}$ the set of all $I\in{\cal I}^k$ such that $wt(I)=w$.
Combining Pieri's formula~(\ref{eq:pieruno})  with Giambelli's formula~(\ref{eq:giambform}), one has, for each $I\in{\cal I}^k$ and each $h\geq 0$:
\[
D_h\Delta_{I}(D)\kformep=D_h\cdot \w^I{\bm\ep}= \sum_{{H\in {\cal P}(I,h)}}
\w^{I+H}{\bm\ep}=\sum_{{H\in {\cal P}(I,h)}}
\Delta_{I+H}(D)\kformep,
\]
proving the equality  
$
D_h\Delta_{I}(D)=\sum_{{H\in {\cal P}(I,h)}}
\Delta_{I+H}(D)
$
in the ring  ${\cal A}^*(\wkM(\ttp))$.

\claim{{\bf Remark} (see~\cite{LakTh1}){\bf.}} \label{rmklakth}  Let  $Split^k_A(\ttp)$ be the universal splitting algebra of the monic polynomial $\ttp$ into the product of two monic polynomials, one of degree $k$. Let  $\ttp={\tt p}_1{\tt q}$ be the universal splitting of $\ttp$ in $Split^k_A(\ttp)$, where $\deg(\ttp_1)=k$, and denote by $s_i$  the complete symmetric polynomial of degree $i$ in the universal roots of $\ttp_1$. Then $Split_A^k(\ttp)$  is generated, as an $A$-algebra, by $(s_i)_{i\geq 1}$ and the map ${\cal A}^*(\wkM(\ttp))\sra Split^k_A(\ttp)$, defined by $D_i\mapsto s_i$, is an $A$-algebra isomorphism. This is because of the module structure of $\bw^kA[X]$ over the ring of symmetric functions defined and studied in~\cite{LakTh}. In fact our formula~(\ref{eq:pieruno}) is the same as Pieri's formula (2.1.1) of~\cite{LakTh}, after replacing $s_i$ with $D_i$.  Let $p:E\sra {\cal Y}$ be a vector bundle of rank $n$ and let $p_k:G(k,E)\sra {\cal Y}$ be the Grassmann bundle over ${\cal Y}$ of $k$-planes in the fibers of $E$. In~\cite{LakTh1} the authors show that, if $A:=A^*({\cal Y})$ is the Chow ring of ${\cal Y}$ and $\ttp=X^n+c_1X^{n-1}+\ldots+c_n\in A[X]$ is such that $c_i:=c_i(E)$ are the Chern classes of $E$,  there is an isomorphism $Split_A^k(\ttp)\sra A^*(G(k,E))$. Let ${\cal Q}_k$ {  be} the universal quotient bundle over $G(k,E)$. Then, the same proof as in~\cite{LakTh1}  works using derivations: by the basis theorem  (\cite{Fu1}, p.~268) the unique $A$-module homomorphism \linebreak $\iota_k: A^*(G(k,E))\sra {\cal A}^*(\wkM(\ttp))$,  mapping $\Delta_I(c({\cal Q}_k-p_k^*E)$ to $ \Delta_I(D)$, is certainly an  isomorphism. To check that it is also a ring homomorphism,  {  it is sufficient} to check it on products of the form $c_h({\cal Q}_k-p_k^*E)\cdot\Delta_I(c({\cal Q}_k-p_k^*E)$: 
\begin{eqnarray*}
\iota_k(c_h({\cal Q}_k-p_k^*E)\cdot \Delta_I(c({\cal Q}_k-p_k^*E))&=&\iota_k\big(\sum_{H\in {\cal P}(I,h)}\Delta_{I+H}(c({\cal Q}_k-p_k^*E))\big)=\\
=\sum_{H\in {\cal P}(I,h)}\Delta_{I+H}(D)=D_h\Delta_I(D)&=&\iota_k(c_h({\cal Q}_k-p_k^*E)\cdot \iota_k(\Delta_I((c({\cal Q}_k-p_k^*E)),
\end{eqnarray*}
by~\cite{Fu1}, Proposition~14.6.1, and~\ref{prel211}.

\claim{\bf Remark.}\label{rmk312} Theorem~\ref{thm1rec} can be proven by showing  that  for each $I\in{\cal I}^k$, there exists $G_{I}\in A[\TT]$ such that
$
\w^I{\bm\ep}=G_{I}(D)\cdot \kformep.
$
This can be achieved via {\em integration by parts} (\ref{eq:intpart2}), as follows.
{  We say} that  $\wkM({\tt p})$ enjoys the property ${\bf G}_j$, for some $1\leq j\leq k$,  if, for each
$
i_{j+1}<\ldots <i_{k}
$
such that $j<i_{j+1}$,
there exists a polynomial $G_{j,i_{j+1},\ldots,i_k}\in A[\TT]$ such that
$
\ep^1\wedge\ldots\wedge \ep^j\wedge \ep^{i_{j+1}}\wedge\ldots\wedge \ep^{i_k}=G_{j,i_{j+1},\ldots,i_k}(D)\cdot \ep^1\wedge\ldots\wedge \ep^k.
$
We shall show, by descending induction, that $\wkM({\tt p})$ enjoys ${\bf G}_j$ for each $1\leq j\leq k$. In fact ${\bf G}_k$ is trivially true. 

Let us suppose that ${\bf G}_j$ holds for some $2\leq j\leq k-1$. Then ${\bf G}_{j-1}$ holds. In fact, for each $j-1<i_{j}<\ldots<i_k$,
\[
\ep^1\wedge\ldots\wedge \ep^{j-1}\wedge \ep^{i_j}\wedge\ldots
\wedge \ep^{i_k}=D_{i_j-j}(\ep^1\wedge\ldots\wedge \ep^{j-1}\wedge \ep^j)\wedge \ep^{i_{j+1}}\wedge\ldots\wedge \ep^{i_k},
\]
basically by~\cite{Gat1},  Corollary~2.5. By {  a}pplying integration by parts~(\ref{eq:intpart2}), one gets:
\[
\ep^1\wedge\ldots\wedge \ep^{j-1}\wedge \ep^{i_j}\wedge\ldots
\wedge \ep^{i_k}=\sum_{h=0}^{i_j-j}D_{i_j-j-h}( \ep^1\wedge\ldots\wedge \ep^j\wedge
\ovD_h( \ep^{i_{j+1}}\wedge\ldots\wedge \ep^{i_k})).
\]
But $\ovD_h( \ep^{i_{j+1}}\wedge\ldots\wedge \ep^{i_k})$ is a sum of elements of the form
$
 \ep^{h_{j+1}}\wedge\ldots\wedge \ep^{h_k}
$, with $j<h_{j+1}<\ldots<h_j$. Then, by the inductive hypothesis, one concludes that ${\bf G}_{j-1}$ holds, too. In particular ${\bf G}_1$ holds and the claim is proven.\qed

\section{Presentations for Intersection Rings}

\begin{claim}{\bf Proposition.}\label{lemma32} {\em Let  ${  D^{-1}_t}:=\displaystyle\sum_{j\geq 0}(-1)^j\ovD_jt^j$ be the inverse of $D_t\in{\cal S}_t(\wM(\ttp))$. 
Then ${\ovD_h}_{|_{\wkM(\ttp)}}=0$, for each $h>k$.
}
\end{claim}
\proof
By induction on $k$.
If $k=0$ one has
$
\ovD_h(m)=0,
$
for each $h\geq 2$ and each $m\in M(\ttp)$. In fact, if $m\in M(\ttp)$, ${D_t}m=\sum_{i\geq 0} D_{1}^im\cdot t^i$. Therefore ${  D_t^{-1}}m=1-D_1m\cdot t$, i.e. ${\ovD_h}_{|_{M(\ttp)}}=0$ for each $h\geq 2$.
Suppose now the property true for $k-1$ and let $h>k$. Any  $m_{k}\in\wkM(\ttp)$ is a finite $A$-linear combination of  elements of the form $m\wedge m_{k-1}$, for suitable $m\in M(\ttp)$ and $m_{k-1}\in\bigwedge^{k-1}M(\ttp)$. It suffices then to check the property for elements of this form. One has:
$
\ovD_h(m\wedge m_{k-1})=\sum_{j=0}^h\ovD_{j}m\wedge \ovD_{h-j}(m_{k-1}).
$
As $\ovD_jm=0$, for $j\geq 2$, it follows that
$
\sum_{j=0}^h\ovD_{j}m\wedge \ovD_{h-j}m_{k-1}=\ovD_1m\wedge \ovD_{h-1}m_{k-1}.
$
By the inductive hypothesis, this last term vanishes as well, because $h-1>k-1$.
\qed

%\claim{} \label{prgiambgen} Recall that  the Schur functions $\{\Delta_\ovI(\TT)\,|\, \ovI \in{\mathbb I}\}$  are a free system of generators of $\ZZ[\TT]$ thought of as a module over the integers. As a consequence $A[\TT]$ is freely generated, as $A$-module, by the Schur functions  $\Delta_\ovI(\TT)$, because $A[\TT]=\ZZ[\TT]\otimes_\ZZ A$.

%

In the sequel $M$ will  be as in~\ref{recalling01} (i.e. $M(\ttp)$ for $\ttp=0$).

\begin{claim}{\bf Proposition.}\label{numgen} {\em The ring ${\cal A}^*(\bw^kM)$ is generated  by $(D_1,D_2,\ldots, D_{k})$ as an $A$-algebra.}
\end{claim}

\proof 
Let ${  D_t^{-1}}=\sum_{i\geq 0}(-1)^i\ovD_it^i$ be the inverse of $D_t$.  First one observes that, for each $h\geq 1$,
$
\ovD_h=\ev_D(\Delta_{(2,3,\ldots,h+1)}(\TT))=\Delta_{(2,3,\ldots,h+1)}(D)
$
and that $\Delta_{(2,3,\ldots,h+1)}(\TT)\in  A[\TT]$ lands in fact in the subring $A[T_1,\ldots, T_h]$ of $A[\TT]$, by Remark~\ref{convx3.8}.
One knows that $\ovD_{k+j}=0$ in ${\cal A}^*(\wkM({\tt p}))$ for each $j\geq 1$ (Proposition~\ref{lemma32}).  Working modulo $\ker(\rho_k)$ (see~\ref{ref4.7}) we may hence write:
\[
\sum_{i\geq 0} D_it^i={1\over 1-\ovD_1t+\ovD_2t^2+\ldots+(-1)^k \ovD_kt^k}.
\]
Define $\widetilde{D}_j(\TT_k)\in A[T_1,\ldots, T_k]\subseteq A[\TT]$ as 
\be
\sum_{j\geq 0} \widetilde{D}_j(\TT_k)t^i={1\over 1-\Delta_{(2)}(\TT)t+\Delta_{(23)}(\TT)t^2+\ldots+ (-1)^k\Delta_{(23\ldots k+1)}(\TT)t^k}.\label{eq:genfunzdtil}
\ee
One clearly has that $D_j-\widetilde{D}_j({\bf D}_k)\in\ker(\rho_k)$ for each $j\geq 0$. Moreover, if $1\leq j\leq k$,
$\widetilde{D}_j(\TT_k)=T_j$, proving  the claim.\qed

\bclm{\bf Example.} In ${\cal A}^*(\Bw^2M)$ one has,
using the recipe~(\ref{eq:genfunzdtil}):

\be
\widetilde{D}_3(\TT_2):=\widetilde{D}_3(T_1,T_2)=T_2T_1-T_1(T_1^2-T_2)=-T_1^3+2T_1T_2\label{eq:d3(d)}
\ee
and
\be
\widetilde{D}_4(\TT_2)=\widetilde{D}_3T_1-T_2(T_1^2-T_2)=(-T_1^3+2T_1T_2)T_1-T_2(T_1^2-T_2)=-T_1^4+T_1^2T_2+T_2^2.\label{eq:d4(d)}
\ee

\eclm
\bclm\label{noreleqnk}{\bf Proposition.} {\em Let  $P\in A[T_1,\ldots,T_k]_w\subset A[\TT]_w$ such that  $P(D)\kformep=0$ ($w\geq 0$). Then $P=0$.
}
\eclm
\proof  
 Any polynomial  $P\in A[T_1,\ldots,T_k]$ of degree $w$ is a unique $A$-linear combination of $\Delta_I(\TT)$, with $I\in{\cal I}^{k,w}$ (since the Schur polynomials $\{\Delta_I(\TT)\,|\, I\in{\cal I}^k\}$ are a $\ZZ$-basis of $\ZZ[\TT]$). Hence
$
P=\sum_{I\in {\cal I}^{k,w}} a_I\Delta_I(\TT)
$
for some (unique!) $a_I\in A_{w-wt(I)}$ and if $P(D)\ep^1\w\ldots\w \ep^{k}=0$, then:
\[
0=P(D)\cdot \ep^1\w\ldots \w \ep^{k}=\sum_{I\in {\cal I}^{k,w}}   a_I\Delta_I(D)\cdot \ep^1\w\ldots \w \ep^{k}=\sum_{I\in {\cal I}^{k,w}} a_I\cdot \w^I{\bm\ep}{  .}
\]
Since $\{\w^I{\bm\ep}\}_{I\in I^{k,w}}$ are $A$-linearly independent, $a_I=0$ for all $I\in {\cal I}^{k,w}$, i.e. $P=0$.
\qed

\bclm\label{presinf}{\bf Corollary.} {\em {  The map~\ref{ref4.7}, $\ev_D:A[\TT]\lra {\cal A}^*(\wM)$ is an isomorphism}. Hence:
\be
{\cal A}^*(\bw M)=A[D]:=A[D_1,D_2,\ldots]\cong A[\TT],\label{eq:presinfinf}
\ee
the polynomial ring in infinitely many indeterminates, while
\be
{\cal A}^*(\bw^kM)=A[{\bf D}_k]:=A[D_1,D_2,\ldots,D_k].\label{eq:presinfk}
\ee
}
\eclm

\proof  Apply~Proposition~\ref{noreleqnk}. One may assume that $P\in A[\TT]$ is homogeneous of degree $w\geq 0$. Suppose that $\ev_D(P)=P(D)=0\in {\cal A}^*(\wM)$.  There is $k\geq 1$ such that $P\in A[T_1,T_2,\ldots, T_{k}]$.  But then
$
P(D)\cdot \kformep=0
$
implies $P=0$, because otherwise one would have a relation (of degree $w$), whence~(\ref{eq:presinfinf}). Since
$
{\cal A}^*(\bw^k M)=\rho_k\big({\cal A}^*(\bw M)\big)
$
and, by Proposition~\ref{lemma32}, $\rho_k(\ovD_h)=0$ for all $h\geq k+1$, one gets the presentation~(\ref{eq:presinfk}).
\qed
%\claim{} Let ${\tt p}\in A[X]$ be a monic polynomial of degree $n$:
%\[
%{\tt p}(X)=X^{n}-c_1X^{n-1}+\ldots+(-1)^{n}c_{n}.
%\]
%Let $M({\tt p})=A[X]/(p)$ and denote by $\ep^{i}=X^{i-1}+(p)$ ($1\leq i\leq i-1$). 
%Then  $M({\tt p})$ is a free $A$-module generated by ${\bm \ep}=(\ep^1,\ldots, \ep^n)$. 

%We shall think of $M({\tt p})$ as a graded $A$-module 
%\[
%M({\tt p})=M({\tt p})_0\oplus M({\tt p})_1\oplus M({\tt p})_2\oplus\ldots=\bigoplus_{h\geq 0}M({\tt p})_h
%\]
%where, for each $h\geq 0$, 
%\[
%M({\tt p})_h= A_h\ep^1\oplus A_{h-1}\ep^2\oplus\ldots\oplus A_0\ep^{h+1}.
%\]
%If $m\in M({\tt p})_h$, one says that $m$ has {\em weight h} and writes $wt(m)=h$; if $a\in A$ is homogeneous, and $m\in M_h$, then $wt(am)=\deg(a)+wt(m)$. In particular $wt(\ep^i)=i-1$. 

% \claim{} Let  ${\bm\nu}=(\nu^1,\nu^2,\ldots)$ defined by
%\[
%\nu^{qn+i}=({\tt p}(X))^qX^{i-1}
%\]
%Then ${\bm\nu}$ is an $A$-basis of $A[X]$ and $\nu^i=X^i$, $1\leq i\leq n$.

%\claim{} One clearly has:
%\[
%\Bw M({\tt p})={\bw A[X]\over \bw A[X]\w {\tt p}(X)A[X]}\qquad and \qquad \Bw^kM({\tt p})={\bw^k A[X]\over \Bw^{k-1}A[X]\w {\tt p}(X)A[X]}
%\]
%where 
%\[
%\bw A[X]\w {\tt p}(X)A[X]=\bigoplus_{k\geq 1}\Bw^{k-1}A[X]\w {\tt p}(X)A[X]
%\] 
%is the bilateral ideal of $\bw A[X]$ generated by ${\tt p}(X)$
%and 
%\[
%\Bw^{k-1}A[X]\w {\tt p}(X)A[X]
%\] 
%is the $A$-submodule of $\bw^kA[X]$ generated by $\ikfornu$, with $i_k>n$.
\claim{}\label{ref5.6} For each $i\geq 1$, let $
\nu^{qn+i}=({\tt p}(X))^qX^{i}
$. Then ${\bm\nu}=(\nu^1,\nu^2,\ldots)$ is an $A$-basis of $M:=M(0)$, such that $\nu^i=X^i$ for each $1\leq i\leq n$.
Let
$
\wM\w {\tt p}M:=\bigoplus_{k\geq 1}\Bw^{k-1}M\w {\tt p}M
$ 
be the bilateral ideal of $\wM$ generated by ${\tt p}$. 
As ${\tt p}M$ is the submodule of $M$ generated by $\nu^i$ with $i>n$, the submodule 
$
\Bw^{k-1}M\w {\tt p}M
$
is the $A$-submodule of $\bw^kM$ generated by $\ikfornu$, with $i_k>n$. The natural map
$
\wM\sra \wM({\tt p})\quad 
$
(resp. $\wkM\sra \wkM({\tt p})$) is surjective and has kernel $\wM\w {\tt p}M$ (resp. $\bw^{k-1}M\w {\tt p}M$).
Hence, one has canonical isomorphisms
\[
\Bw M({\tt p})={\wM\over \wM \w {\tt p}M}\qquad and \qquad \wkM({\tt p})={\wkM\over \Bw^{k-1}M\w {\tt p}M}.
\]
Let $\phi_k:\wkM\sra\wkM(\ttp)$ be  the canonical projection and  let 
\[
J_k({\tt  p}):=\{ P(D)\in A[D_1,\ldots, D_k]\,|\, P(D)\ep^1\w\ldots \w \ep^{k}\in \bw^{k-1}M\w {\tt p}M\},
\]
{  which} is an ideal of $A[D_1,\ldots, D_k]={\cal A}^*(\wkM)$.

\bclm{\bf Theorem.}\label{thme1310} {\em For each $j\geq 1$, let
\[
\widetilde{D}_{n-k+j}({\bf D}_k,\ttp)=\widetilde{D}_{n-k+j}({\bf D}_k)+\sum_{i=1}^{n-k+j}c_i\widetilde{D}_{n-k+j-i}({\bf D}_k).
\]
Then:
\be
J_k({\tt p})=(\widetilde{D}_{n-k+1}({\bf D}_k,\ttp),\ldots,\widetilde{D}_{n}({\bf D}_k,\ttp)).
\label{eq:preintring}
\ee
}
\eclm
\proof
%Notation as in Section~\ref{setup41}. 
Let $D'_t=\sum_{i\geq 0}D'_it^i$ be the unique derivation on $\wM$  such that $D'_t\nu^j=\sum_{i\geq 0}\nu^{i+j}t^i$. Then 
\begin{enumerate}
\item[i)] $D'_i\in {\cal A}^*(\wM(\ttp))$ for each $i\geq 0$, 
\item[ii)] $\rho_k(D'_i)=D_i$ if $1\leq i\leq n-k$ and 
\end{enumerate}
\be
{\rm iii)} \hskip5pt\rho_k({D}'_{n-k+j})=\widetilde{D}_{n-k+j}({\bf D}_k,\ttp),\qquad \forall j\geq 1.\hskip141pt \label{eq:eqpres}
\ee
To check i), is sufficient to show that each $D'_i$ is an $A$-polynomial expression in the $D_i$s. As a matter of fact,  if $i\leq n-k$:
\begin{eqnarray}
D'_i(\kformX)&=&D'_i(\nu^1\w\ldots\w\nu^k)=\nu^1\w\ldots\w\nu^{k-1}\w\nu^{k+j}=\nonumber\\
&=&X^1\w\ldots\w X^{k-1}\w X^{k+i}=D_i(\kformX),\label{eq:final1}
\end{eqnarray}
and,  for each $j\geq 1$:
\begin{center}
$
{D}'_{n-k+j}(\kformX)={D}'_{n-k+j}\nu^1\w\ldots\w\nu^k=\nu^1\w\ldots\w\nu^{k-1}\w \nu^{n+j}=
$

\medskip
$
=X^1\w\ldots\w X^{k-1}\w (X^{n+j}+c_1X^{n+j-1}+\ldots+c_{n+j-k}X^{k}+\ldots+c_{n}X^{j})=
$

\medskip
$
=X^1\w\ldots\w X^{k-1}\w X^{n+j}+\sum_{i=1}^{n-k+j} c_i(X^1\w\ldots\w X^{k-1}\w X^{n+j-i})=
$
\be
{\textstyle =(\tilde{D}_{n-k+j}({\bf D}_k)+\sum_{i=1}^{n-k+j}c_i\tilde{D}_{n-k+j-i}({\bf D}_k))\kformX.}\label{eq:final2}
\ee
\end{center}
Therefore formulas~(\ref{eq:final1}) and~(\ref{eq:final2}) show i), (\ref{eq:final1}) shows ii) and~(\ref{eq:final2}) shows~iii) above.
We can now prove equality~(\ref{eq:preintring}). 
Clearly  $(\widetilde{D}_{n-k+j}({\bf D}_k,\ttp))_{j\geq 1}\subseteq J_k({\tt p})$,
because:
\begin{eqnarray*}
\widetilde{D}_{n-k+j}({\bf D}_k,\ttp)\kfornu&=&D'_{n-k+j}\kfornu=\\&=&\nu^1\w\ldots\w\nu^{k-1}\w\nu^{n+j}\in \Bw^{k-1}M\w {\tt p}M.
\end{eqnarray*}
To show that  $J_k({\tt p})\subseteq (\widetilde{D}_{n-k+j}({\bf D}_k,\ttp))_{j\geq 1}$ as well,  let $P\in A[T_1,\ldots,T_k]\subseteq A[\TT]$ such that $P({ D}')\kformX\in \Bw^{k-1}M\w {\tt p}M$. Without loss of generality one may assume that $P$ is homogeneous of degree $w$.
Then
\[
P(D')\kfornu=\sum a_I\Delta_I(D')\kfornu=\sum \ikfornu,
\]
where last sum is over all $(i_1,\ldots, i_k)\in{\cal I}^{k,{  w}}$ such that $i_k>n$.  By~\ref{convx3.8},   $\Delta_I(D')$  belongs to the ideal $({D}'_{i_k-1},\ldots,{D}'_{i_k-k})$ and, since $i_k>n$, one sees that  if $\Delta_I(D')\kfornu \in \Bw^{k-1}M\w {\tt p}M$, then $\Delta_I(D')\in (\widetilde{D}_{n-k+j}({\bf D}'_k))_{j\geq 1}$.  The relation
\[
D'_{n+1}-D'_n\ovD'_1+\ldots+(-1)^{n-k+1} D'_{n-k+1}\ovD'_k=0,
\]
 holding in ${\cal A}^*(\wkM(\ttp))$, implies that ${D}'_{n+1}\big(\wkM(\ttp)\big)\in ({D}'_{n-k+1},\ldots,{D}'_{n})\wkM(\ttp)$ (here $(-1)^i\ovD'_i$, as in~\ref{conv3.9},  stands for the $i^{th}$ coefficient of ${(D_t')}^{-1}$). 
 
 By induction ${D}'_{n+j}\big(\wkM(\ttp)\big)\in({D}'_{n-k+1},\ldots,{D}'_{n})\bw^kM$ as well. Because of~(\ref{eq:eqpres}), one hence has 
$
J_k(\ttp)\subseteq (\rho_k(D'_{n-k+1}),\ldots,\rho_k(D'_{n}))=(\widetilde{D}_{n-k+1}({\bf D}_k,\ttp),\ldots, \widetilde{D}_{n}({\bf D}_k,\ttp)),
$
i.e. $J_k({\tt p})$ is given precisely by~(\ref{eq:preintring}).\qed

\bclm{\bf Theorem.} \label{generalpres}{\em The following isomorphism holds:
\be
A^*(\wkM({\tt p}))={A[D_1,\ldots, D_k]\over (\widetilde{D}_{n-k+1}({\bf D}_k,\ttp),\ldots, \widetilde{D}_{n}({\bf D}_k,\ttp))}.\label{eq:presfinmp}
\ee
}
\eclm
\proof  Notation as in~\ref{ref5.6}. {  Recall that by Corollary~\ref{presinf}, formula~(\ref{eq:presinfk}), $D_1,\ldots, D_k$, are algebraically independent elements of ${\cal A}^*(\wkM)$}. Clearly, $P(D)\in\ker(\phi_k)$ if and only if $P(D)\kformep\in \bw^{k-1}M\w \ttp M$, i.e. if and only if $P(D)\in J_k(\ttp)$. Hence 
$
{\cal A}^*(\wkM(\ttp))={{\cal A}^*(\wkM)/J_k(\ttp)},
$
 and the conclusion follows by Corollary~\ref{presinf} and Theorem~\ref{thme1310}.\qed

\claim{\bf Remark.} By~\cite{LakTh1},  the polynomial $\ttp$, as above, splits in the ring ${\cal A}^*(\wkM(\ttp))[X]$  as the product  $(X^k+D_1X^{k-1}+\ldots+D_k)\cdot {\tt q}$, where ${\tt q}$ is a monic polynomial of degree $n-k$ (with ${\cal A}^*(\wkM(\ttp))[X]$-coefficients).

\claim{\bf Examples.} \label{ex89}
Let $M=XA[X]$ and ${\tt p}(X)=X^4+c_1X^3+c_2X^2+c_3X+c_4$, $c_i\in A_i$.
Then one has:
\be
{\cal A}^*(\Bw^2M({\tt p}))={A[D_1,D_2]\over(\widetilde{D}_3({\bf D}_2,\ttp),\widetilde{D}_4({\bf D}_2,\ttp))},\label{eq:irg24}
\ee
where 

\smallskip
 \noindent $\widetilde{D}_3({\bf D}_2,\ttp)=\widetilde{D}_3({\bf D}_2)+c_1\widetilde{D}_2({\bf D}_2)+c_2\widetilde{D}_1({\bf D}_2)+c_3=2D_1D_2-D_1^3+c_1D_2+c_2D_1+c_3$

\smallskip
\noindent
and
\begin{center}
\smallskip
\noindent
$
\widetilde{D}_4({\bf D}_2,\ttp)=\widetilde{D}_4({\bf D}_2)+c_1\widetilde{D}_3({\bf D}_2)+c_2\widetilde{D}_2({\bf D}_2)+c_3\widetilde{D}_1({\bf D}_2)+c_4=$

\smallskip
\hskip68pt $=D_2^2+D_1^2D_2-D_1^4+c_1(2D_1D_2-D_1^3)+c_2D_2+c_3D_1+c_4$,
\end{center}
which we obtained from~(\ref{eq:d3(d)}) and~(\ref{eq:d4(d)}). Let us enumerate some particular cases.

\item[1)] If $A=\ZZ$, thought of as a graded ring concentrated in degree $0$, then $c_i=0$, $1\leq i\leq 4$. Then ${\tt p}=X^4$ and presentation~(\ref{eq:irg24}) becomes:
\[
{\cal A}^*(\bw^2M(X^4))={A[D_1,D_2]\over(2D_1D_2-D_1^3, D_2^2+D_1^2D_2-D_1^4)},
\]
which coincides (Cf.~\cite{Gat1},~\cite{Gat2}) with  the presentation of the integral cohomology ring of the grassmannian $G(2,4)$ of $2$-planes in $\CC^4$ (or of the grassmannian of lines $G(1,\PP^3)$ in the complex projective $3$-space).
\item[2)]
If $A=\ZZ[q]$,  and ${\tt p}(X)=X^4+q$, then~(\ref{eq:irg24}) reads:
\[
{\cal A}^*(\bw^2M(X^4+q))={\ZZ[q][D_1,D_2]\over(2D_1D_2-D_1^3, D_2^2+D_1^2D_2-D_1^4+q)},
\]
which is the Witten-Siebert-Tian presentation of the {\em small quantum cohomology ring} $QH^*(G(2,4))$ (\cite{Wi},~\cite{SiebTi},~\cite{Ber1}; see also~\cite{Gat1});
\item[3)]  If $\pi:E\sra {\cal Y}$ is a holomorphic vector bundle of rank $4$ on a smooth complex variety of dimension $m\geq 0$, and ${\tt p}(X)= X^4+\pi^*c_1X^3+\pi^*c_2X^2+\pi^*c_3X+\pi^*c_4\in A^*({\cal Y})[X]$, where $c_i$ are the Chern classes of $E$ as in~\cite{Fu1}, p.~141, then $M({\tt p})=A_*(\PP(E))$, $A=A^*({\cal Y})$ and $D_1=c_1(O_{\PP(E)}(-1))$, thought of as operator on $A_*({\cal Y})$;   in this case~(\ref{eq:irg24}) gives the presentation of  $A^*(G(2,E))$ (Cf.~\ref{rmklakth}). If ${\cal Y}$ is a point, then $A^*({\cal Y})=\ZZ$, $c_i=0$ and one recovers once again the presentation of the Chow ring of the grassmannian $G(2,4)$.
\item[4)] Let $A=\ZZ[y_1,y_2,y_3,y_4]$ and 
$
\ttp(X)=\prod_{i=1}^4(X-y_i+y_1)+q\in A[q].
$
In this case presentation~(\ref{eq:irg24}) is that the {\em quantum equivariant cohomology ring} $QH^*_T(G(2,4))$  of the Grassmannian $G(2,4)$ under the action of a $4$-dimensional compact or algebraic torus via a diagonal action with only isolated fixed points, as studied by Mihalcea in~\cite{Mihalcea2}, Theorem 4.2, setting $p=2$ and  $m=4$. This is compatible with the main result of the paper~\cite{GatSant2} (Theorem~3.7), with Theorem~2.9 of~\cite{Gat1} and is now a consequence of~\cite{formalism}. Notice that our generators are not the same as used in~\cite{Mihalcea2} (Cf.~\cite{GatSant2}).

%\newpage

%\input bibliografiaSCGA

\medskip

Dipartimento di Matematica, Politecnico di Torino, C.so Duca degli Abruzzi, 24, 

10129, Torino

\medskip
{\tt letterio.gatto@polito.it}
\hfill{\tt taise@calvino.polito.it}

\end{document}